\documentclass[11pt]{article} \usepackage[dvips]{epsfig}
\usepackage{latexsym}
\usepackage{amsfonts,amsmath,amssymb}

\newtheorem{lemma}{Lemma}

\newtheorem{theorem}{Theorem}
\newtheorem{proposition}{Proposition}
\newtheorem{question}{Question}
\newtheorem{definition}{Definition}
\newcommand{\qed}{\mbox{$\Diamond$}\vspace{\baselineskip}}
\newenvironment{proof}{\noindent{\bf Proof:}}{\qed}

\begin{document}

\title{Where the monotone pattern (mostly) rules}

\author{Mikl\'os B\'ona \\
Department of Mathematics \\
University of Florida\\
 Gainesville FL 32611-8105\\
bona@math.ufl.edu \thanks{Partially supported
by an NSA Young Investigator Award.}}

\date{}

\maketitle

\begin{abstract} We consider pattern containment and
avoidance with  a very
tight definition that was used first by Riordan more than 60 years ago.
Using this definition, we prove the monotone pattern is easier to avoid
than almost any other pattern of the same length.
 We also show that with this definition, almost
all patterns of length $k$ are avoided by the same number of permutations
of length $n$. The corresponding statements are not known
to be true for more relaxed definitions of pattern containment. This is the
first time we know of that expectations are used to compare numbers of
permutations avoiding certain patterns. 
\end{abstract}

\section{Introduction}
The classic definition of pattern avoidance on permutations is as follows.
Let $p=p_1p_2\cdots p_n$ be a permutation, let $k<n$, and let $q=q_1q_2\cdots
q_k$ be another permutation. We say that $p$ {\em contains} $q$ as a pattern
if there exists a subsequence $1\leq i_1<i_2<\cdots <i_k\leq n$
 so that for all indices
$j$ and $r$, the inequality $q_j<q_r$ holds if and only if the inequality
$p_{i_j}<p_{i_r}$ holds. If $p$ does not contain $q$, then we say
that $p$ {\em avoids} $q$.  In other words, $p$ contains $q$ if $p$ has
a subsequence of entries, not necessarily in consecutive positions, which
relate to each other the same way as the entries of $q$ do. 

Classic pattern avoidance has been a rapidly developing field for the last
decade. One of the most fascinating subjects in this field was the enumeration
of permutations avoiding a given pattern. Let $S_n(q)$ denote the number
of permutations of length $n$ (or in what follows, $n$-permutations)
 that avoid the pattern $q$, and let us
consider the numbers $S_n(q)$ for each pattern $q$ of length $k$. A very 
interesting and counter-intuitive phenomenon is that in this multiset of
$k!$ numbers, the number $S_n(q)$ will, in general, not be the largest or
the smallest number. There are several results on this fact (see 
\cite{bona}, \cite{bonal} or \cite{records}), but the 
phenomenon is still not perfectly well understood. 

In 2001, Elizalde and Noy \cite{elizalde}
 proposed another definition of pattern containment.
We will say that the permutation $p=p_1p_2\cdots p_n$ {\em tightly} 
contains the permutation  $q=q_1q_2\cdots q_k$ if there exists an index
$0\leq i\leq n-k$ so that $q_j<q_r$ if and only if $p_{i+j}<p_{i+r}$. 
In other words, for $p$ to contain $q$, we require that $p$ has a {\em
consecutive} string of entries that relate to each other the same way the 
entries of $q$ do. For instance, 246351 contains 132 (take the second, third,
and fifth entries, for instance), but it does not {\em tightly contain}
132 since there are no three entries in consecutive positions in 246351
that would form a 132-pattern.
 If $p$ does not tightly contain $q$, then we say that $p$ {\em
tightly avoids} $q$. Let $T_n(q)$ denote the number of $n$-permutations
that tightly avoid $q$. Elizalde and Noy conjectured in \cite{elizalde}
 that  no pattern of
length $k$ is tightly avoided by more $n$-permutations than the monotone
pattern. In other words, if $q$ is a pattern of length $k$, then
\begin{equation} \label{elinoy} T_n(q)\leq T_n(12\cdots k).\end{equation}
This conjecture is still open. (In the special case of $k=3$,
it was proved in \cite{elizalde}.) Still, it is worth pointing out
that changing the definition of pattern avoidance changed the status of the
monotone pattern among all patterns of the same length. With this definition,
it is believed that the monotone pattern is the easiest pattern to avoid.    

This perceived change in the status of the monotone pattern led us to the
following direction of research. Let us take the idea of Elizalde and Noy
one step further, by restricting the notion of pattern containment further
as follows. Let $p=p_1p_2\cdots p_n$ be a permutation, let $k<n$, and let
 $q=q_1q_2\cdots
q_k$ be another permutation. We say that $p$ {\em very tightly} contains
$q$ if there is an index $0\leq i\leq n-k$ and an integer $0\leq a\leq 
n-k$ so that  $q_j<q_r$ if and only if $p_{i+j}<p_{i+r}$, and, 
\[\{p_{i+1},p_{i+2},\cdots ,p_{i+k}\}=\{a+1,a+2,\cdots ,a+k\} .\]
That is, $p$ very tightly contains $q$ if $p$ tightly contains $q$ and
the entries of $p$ that form a copy of $q$ are not just in consecutive
positions, but they are also consecutive as integers (in the sense that 
their set is an interval).

For example, 15324 tightly contains 132 (consider the first three entries),
but does not very tightly contain 132. On the other hand, 15324 very tightly
contains 213, as can be seen by considering the last three entries. If $p$
does not very tightly contain $q$, then we will say that $p$ {\em very tightly
avoids} $q$. Note that in the special case when $q$ is the monotone pattern,
this notion was  studied before pattern avoidance became widely known.
The literature of permutations very tightly avoiding monotone patterns
goes back at least to \cite{Riordan}. More recent
 examples include \cite{Jackson} and \cite{JackReid}. However, we did not find
any examples where the notion was used in connection with any other pattern.

Let $V_n(q)$ denote the number of $n$-permutations that very tightly avoid
$q$. While we cannot prove that $V_n(q)\leq V_n(12\cdots k)$ for all
patterns $q$ of length $k$, we will be able to prove that this inequality
holds for {\em most} patterns $q$ of length $k$. As a byproduct, 
we will prove that for all $k$, 
there exists a set $W_k$ of patterns of length $k$ so 
that $\lim_{k\rightarrow \infty}\frac{|W_k|}{k!}=1$, and 
 $V_n(q)$ is identical for all patterns $q\in W_k(q)$. In other words, almost
all patterns of length $k$ are equally difficult to very tightly
 avoid. There are no 
comparable statements known for the other two discussed notions of pattern
avoidance.

Our argument will be a probabilistic one. Once the framework is set up, 
the computation will be elementary. However, this is the first time
we know of that expectations are successfully 
used to compare the number of permutations avoiding a given pattern 
(admittedly, with a very restrictive definition of pattern avoidance). We
wonder whether more sophisticated methods of enumeration could extend 
the reach of this technique to less restrictive definitions of pattern
avoidance.

\section{A Probabilistic Argument}
\subsection{The outline of the argument}
For the rest of this section, let $k\geq 3$ be a fixed positive integer.
Let $\alpha=12\cdots k$, the monotone pattern of length $k$.
Recall that  $V_n(\alpha)$ is
 the number of $n$-permutations very tightly
 avoiding $\alpha$. Our goal is true prove that
\[V_n(q)\leq V_n(\alpha)\] for any pattern $q$ of length $k$.

Let $q$ be any pattern of length $k$.
For a fixed positive integer  $n$, let $X_{n,q}$ be the random 
variable counting the occurrences
of $q$ in a randomly selected 
$n$-permutation. As the following 
straightforward proposition shows, the expectation of $X_{n,q}$ does not
depend on $q$; it only depends on $n$, and the length $k$ of $q$.

\begin{proposition} \label{equal}
For any fixed $n$, and $q\in S_k$,  we have
\[E(X_{n,q})=\frac{(n-k+1)^2}{{n\choose k}k!}.\]
\end{proposition}

\begin{proof} Let $X_i$ be the indicator random variable of the event that
the string $p_{i+1}\cdots p_{i+k}$ is a $q$-pattern in the very tight sense.
 Then 
$E(X_i)=P(p_{i+1}\cdots p_{i+k} \simeq q)=\frac{n-k+1}{{n\choose k}}
\cdot \frac{1}{k!}$, since there are $n-k+1$ favorable choices for the set
of the entries $p_{i+1},\cdots, p_{i+k}$, and there is $1/k!$ chance that
their pattern is $q$. Now note that $E(X _{n,q})=\sum_{i=0}^{n-k}
X_i$, and the statement is proved by the linearity of expectation.
\end{proof}

Let $p_{n,i,q}$
 be the probability that a randomly selected $n$-permutation
contains {\em exactly} $i$ copies of $q$, and let $P_{n,i,q}$ be the
probability that a randomly selected $n$-permutation contains {\em at least} 
$i$  copies of $q$.

Set $m=n-k+1$, and observe
 that no $n$-permutation can very tightly contain more
than $m$ copies of any given pattern $q$ of length $k$. 
Now note that by the definition of expectation
\begin{eqnarray*} \label{atleast} 
E(X_n,q) & = & \sum_{i= 1}^{m} ip_{n,i,q} \\
  & = & \sum_{j=0}^{m-1} \sum_{i=0}^j p_{n,m-i,q} \\
 & = & p_{n,m,q} + (p_{n,m,q}+p_{n,m-1,q}) + \cdots + 
(p_{n,m,q}+\cdots +p_{n,1,q}) \\
 & = & \sum_{i=1}^m P(n,i,q) .
\end{eqnarray*}

By Proposition $\ref{equal}$, we know that $E(X_{n,q})=E(X_{n,\alpha})$,
and then previous displayed equation implies that
\begin{equation}\label{bigequal}
\sum_{i=1}^m P(n,i,q) = \sum_{i=1}^m P(n,i,\alpha) .\end{equation}
So if we could show that for $i\geq 2$, the inequality 
\begin{equation} \label{toprove} P(n,i,q)\leq P(n,i,\alpha)
\end{equation} holds, then (\ref{bigequal}) would imply
that  $P(n,1,q) \geq P(n,1,\alpha)$, which is just what we set out
to prove.  

The simple counting argument that we present in this paper will not prove
 (\ref{bigequal}) for every pattern $q$.
However, it will prove  (\ref{bigequal}) 
 for most patterns $q$. We describe these
patterns in the next subsection.

\subsection{Condensible and Non-condensible Patterns}
Let us assume that the permutation $p=p_1p_2\cdots p_n$ very tightly contains
two {\em non-disjoint} copies of the pattern $q=q_1q_2\cdots q_k$.
Let these two copies be  $q^{(1)}$ and
 $q^{(2)}$, so that $q^{(1)}=p_{i+1}p_{i+2}\cdots p_{i+k}$ and 
$q^{(2)}=p_{i+j+1}p_{i+j+2}\cdots p_{i+j+k}$ for some $j\in [1,k-1]$.
 Then $|q^{(1)}\cap q^{(2)}|=k-j+1=s$. Furthermore, since the set of entries 
of $q^{(1)}$ is an interval, and the set of entries of $q^{(2)}$ is  an
interval, it follows that the set of entries of $q^{(1)}\cap q^{(2)}$ is
also an interval. So the rightmost
 $s$ entries of $q$, and the leftmost $s$ entries
of $q$ must form identical patterns, and the respective sets of these entries
must both be intervals. 

For obvious symmetry reasons, we can assume that $q_1<q_k$. We claim that then
the {\em rightmost} $s$ entries of $q$ must also be the {\em largest} $s$ 
entries of $q$. This can be seen by considering  $q^{(1)}$. Indeed,
 the set of these entries of   $q^{(1)}$ is the 
intersection of two intervals of the same length, and therefore, must
 be an ending segment of
the interval that starts on the left of the other. An analogous argument, 
applied for  $q^{(2)}$, shows that the leftmost $s$ entries of $q$ must
also be the {\em smallest} $s$ entries of $q$. 

The following Proposition collects the observations made in this subsection.

\begin{proposition} \label{conditions}
Let $p$ be a permutation that very tightly  contains copies $q^{(1)}$ and
 $q^{(2)}$ of the pattern $q=q_1q_2\cdots q_k$. Let us assume that
$q_1<q_k$. Then 
 $q^{(1)}$ and $q^{(2)}$ are disjoint unless all of the following hold.

There exists a positive integer $s\leq k-1$ so that 
\begin{enumerate}
\item the rightmost $s$ entries of $q$ are also the largest $s$ entries 
of $q$,
and the leftmost $s$ entries of $q$ are also the smallest $s$ entries of $q$,
and 
\item the pattern of the leftmost $s$ entries of $q$ is identical to the 
pattern of the rightmost $s$ entries of $q$. 
\end{enumerate}
\end{proposition}

It is easy to see that if $q$ satisfies both of these criteria, then two
very tightly contained copies of $q$ in $p$ may indeed intersect. 
For example, the pattern $q=2143$ satisfies both of the above criteria with
$s=2$, and indeed, 214365 very tightly contains two intersecting copies of
$q$, namely 2143 and 4365. 

\begin{definition} 
Let $q$ be a pattern that satisfies both conditions of Proposition 
\ref{conditions}. Then we say that $q$ is condensible. 
\end{definition}

It is not difficult to prove that almost all patterns of length $k$ are 
non-condensible. We do not want to break the course of our proof with 
this, so we postpone this computation until the Appendix. 

\subsection{The Computational Part of the Proof}
The following Lemma is the heart of our main result.
\begin{lemma} \label{heart} Let $q$ be a non-condensible pattern, and let 
$i>1$. Then \[P(n,i,q)\leq P(n,i,\alpha).\]
\end{lemma}

\begin{proof}
We point out that if $n<ik$, then the statement is clearly true. Indeed, 
$P(n,i,q)=0$ since any two copies of $q$ in an $n$-permutation $p$
would have to be disjoint, and $n$ is too small for that. Therefore, in the
rest of the proof, we can 
assume that $n\geq ik$. For $k\leq 2$, the statement is
trivial, so we assume $k\geq 3$ as well.  

First, we prove a lower bound on $P(n,i,\alpha)$.
The number of $n$-permutations very tightly
containing $i$ copies of
$\alpha$ is at least as large as the number of $n$-permutations very tightly
containing
the pattern $12\cdots (i+k-1)$. The latter is at least as large as the number
of $n$-permutations that very tightly contain
 a $12\cdots (i+k-1)$-pattern in their first $i+k-1$ positions, that is,
$(n-k-i+2)\cdot (n-k-i+1)!=(n-k-i+2)!$. Therefore, 
\begin{equation} \label{lower}
 \frac{(n-k-i+2)!}{n!} \leq P(n,i,\alpha).
\end{equation}

We are now going to find an upper bound for $P(n,i,q)$.
 Let $S$ be an $i$-element
subset of $[n]$ so that the elements of $S$ can be the starting positions
of $i$ (necessarily disjoint)
 very tight copies of $q$ in an $n$-permutation. If
$S=\{s_1,s_2,\cdots ,s_i\}$, then 
this is equivalent to saying that 
\[1\leq s_1< s_2-k+1 \leq s_3-2k+2\leq \cdots \leq s_i-(i-1)(k-1)\leq n-
i(k-1).\]
Therefore, there are ${n-i(k-1)\choose i}$ possibilities for $S$. 
Now let $A_S$ be the event that in a random permutation $p=p_1\cdots p_n$,
the subsequence $p_jp_{j+1}\cdots p_{j+k-1}$ is a very tight
 $q$-subsequence for
all $j\in S$. Let $A_{i,q}$ be the event that $p$ contains at
least $i$ very tight copies of $q$. Then $P( A_{i,q})=P(n,i,q)$. 
Furthermore, 
\[A_{i,q}= \cup_S A_S,\]
where the union is taken over all  ${n-i(k-1)\choose i}$ possible subsets
for $S$. Therefore, 
\begin{equation}
\label{estimate} P(n,i,q)=P( A_{i,q})\leq \sum_S P(A_S).
\end{equation}
Let us now compute $P(A_S)$. We will see that this probability does not
depend on the choice of $S$. Indeed, just as there are ${n-i(k-1)\choose i}$
 possibilities for $S$, there are ${n-i(k-1)\choose i}$ possibilities for 
the {\em entries} in the positions belonging to $S$. Once those entries
are known, the rest of the $q$-patterns starting in those entries
are determined, and there are $(n-ik)!$ possibilities for the rest
of the permutation. This shows that $P(A_S)={n-i(k-1)\choose i}(n-ik)!
\frac{1}{n!}$ for
all $S$. Therefore, (\ref{estimate}) implies
\begin{equation} \label{upper} 
P(n,i,q)\leq {n-i(k-1)\choose i}^2(n-ik)!\frac{1}{n!}.\end{equation}

Comparing (\ref{upper}) and (\ref{lower}), we see that
 our lemma will be proved if we show that for $i>1$, the inequality
\[{n-i(k-1)\choose i}^2(n-ik)! \leq (n-i-k+2)!,\]
or, equivalently, 
\begin{equation}
\label{inequal} (n-i(k-1))_i\leq i!^2 (n-k-i+2)(n-k-i+1)\cdots (n-i(k-1)+1))
\end{equation} holds. 
Where $(z)_j=z(z-1)\cdots (z-j+1)$. 
Note that the left-hand side has $i$ factors, while the right-hand side, not
counting $i!^2$, has $(k-1)(i-1)>i$ factors, each of which are larger
than the factors of the left-hand side. Therefore, (\ref{inequal}) holds,
and the Lemma is proved. 
\end{proof}

The proof of our main result is now immediate.

\begin{theorem} Let $q$ be any pattern of length $k$. Then
\[V_n(q)\leq V_n(\alpha).\]
\end{theorem}

\begin{proof} Lemma \ref{heart} and formula (\ref{bigequal}) together imply
that $P(n,1,q)\geq P(n,1,\alpha)$, which means that there are 
at least as many 
$n$-permutations that very tightly contain $q$ as $n$-permutations that
very tightly contain $\alpha$. 
\end{proof}

\subsection{A Result on Non-condensible Patterns}
We have seen in Proposition \ref{conditions} that if $q$ is non-condensible
and $q_1<q_k$,
then any two copies of $q$ contained in a given permutation $p$ are disjoint.
Therefore, the number of $n$-permutations that very tightly avoid $q$ can
be computed by the Principle of Inclusion-Exclusion. 
Indeed, in this case, the following holds.

\begin{proposition} Let $q$ be a non-condensible pattern. Then
\[V_n(q)=n! - \sum_{i=1}^{\lfloor n/k \rfloor} {n-i(k-1)\choose i}^2
(n-ik)!.\] In particular, $V_n(q)$ does not depend on the choice of 
$q$.  
\end{proposition}

\begin{proof} 
In the proof of Lemma \ref{heart}, more precisely, in our argument showing
that (\ref{upper}) holds, we showed that there are ${n-i(k-1)\choose i}$ 
ways to choose an $i$-element set of positions that can be the starting
positions of $i$ disjoint very tight copies of $q$, and there are 
 ${n-i(k-1)\choose i}$ ways to choose the sets of entries forming these
same copies. Once these choices are made, the rest of the permutation
can be chosen in $(n-ik)!$ ways. The statement now follows by the 
Principle of Inclusion-Exclusion.
\end{proof}

As we said, we will prove in Proposition \ref{almostall} that almost
all patterns are non-condensible.
Note that nothing comparable is known for the other two notions of
pattern avoidance. Numerical evidence suggests that similarly strong
results will probably {\em not} hold if the traditional definition or the
tight definition is used. 

\section{Further Directions}

The novelty of this paper, beside the notion of very tight containment,
was the application of expectations to compare the numbers of permutations
avoiding two patterns. The computations themselves were elementary. 
This leads to the following question.

\begin{question} Is it possible to apply our method to compare the numbers
$T_n(q)$ and $T_n(\alpha)$, or the numbers $S_n(q)$ and $S_n(\alpha)$, for at
least some patterns $q$?
\end{question}

As we used a very simple estimate in our proof of Lemma \ref{heart}, 
there may be room for improvement at that point. 

There are several natural  questions that can be raised about the enumeration
of permutations that very tightly avoid a pattern. Let us recall that we
proved a formula for $V_n(q)$ for the overwhelming majority of patterns
$q$, namely for non-condensible patterns. A formula for $V_n(\alpha)$, 
where $\alpha$ is the monotone pattern, can be found in \cite{Jackson} and
\cite{JackReid}. This raises the following question. 

\begin{question}
Are there other patterns $q$ for which $V_n(q)$ can be explicitly determined?
\end{question}

Let us call patterns $q$ and $q'$ {\em very tightly equivalent} if 
$V_n(q)=V_n(q')$ for all $n$. We have seen that almost all patterns of
length $k$ are very tightly equivalent. This raises the following questions.

\begin{question}
How many equivalence classes are there for very tight patterns of length $k$?
\end{question}

\begin{question} Can we say anything about the same topic for tight pattern
containment, or traditional pattern containment? If equivalence is defined
for them in an analogous way, how many equivalence classes will be formed, 
(for patterns of length $k$)
and how large will the largest one be?
\end{question}

\section{Appendix}

In this section we prove the following simple fact.

\begin{proposition} \label{almostall}
Let $h_n$ be the number of condensible permutations of length $n$.
Then \[\lim_{n\rightarrow \infty} \frac{h_n}{n!}=0.\]
\end{proposition}

\begin{proof}
We prove that another class of permutations, one that contains all
condensible permutations, is also very small. 
Let 
 $a_n$ be the number of  {\em decomposable} permutations, that is,
permutations $p_1p_2\cdots p_n$ for which is there is an index $i$ so that
$p_j<p_m$ if $j\leq i<m$. In other words, $p$ can  be cut into two parts
so that everything before the cut is less than everything after the cut. 
(Note that if a permutation is not decomposable, then it is called 
{\em indecomposable}, and Exercise 1.32 of \cite{stanley} contains more
information about these permutations.)

Counting according to the index $i$ of the above definition, we see that
\begin{eqnarray*} \frac{a_n}{n!} & \leq & \sum_{i=1}^{n-1}
 \frac{i!(n-i)!}{n!} \\ 
& \leq & 
\sum_{i=1}^{n-1} \frac{{n\choose i}}{n} \\
& \leq & \frac{2}{n} + \frac{n-3}{{n\choose 2}},
\end{eqnarray*}
where in the last step we used the well-known unimodal property of binomial
coefficient, in particular the inequality that ${n\choose i}\geq {n\choose 2}$
if $2\leq i\leq n-2$. 

Therefore, $a_n/n!$ converges to 0 as $n$ goes to infinity. Clearly, all
condensible patterns are decomposable since their first $i$ entries are
also their $i$ smallest entries, so our claim  follows. 
\end{proof}

\end{document}